\documentclass[a4paper,11pt,reqno,noindent]{amsart}
\usepackage[centertags]{amsmath}
\usepackage{amsfonts,amssymb,amsthm,dsfont,cases,amscd,esint,enumerate}
\usepackage[T1]{fontenc}
\usepackage[english]{babel}
\usepackage[applemac]{inputenc}
\usepackage{newlfont}
\usepackage{color}
\usepackage[body={15cm,21.5cm},centering]{geometry} 
\usepackage{fancyhdr}
\usepackage{enumerate}

\usepackage[colorlinks,citecolor=black,urlcolor=black]{hyperref}

\theoremstyle{plain}
\newtheorem{theor0}{Theorem}[section]

\newtheorem{lem0}[theor0]{Lemma}

\newtheorem{prop0}[theor0]{Proposition}
\newenvironment{prop}
  {\pushQED{\qed}\begin{prop0}}{\popQED\end{prop0}}
\newtheorem{cor0}[theor0]{Corollary}

\newtheorem{propr0}[theor0]{Property}

\newtheorem{hyp0}[theor0]{Hypothesis}

\newtheorem{result0}[theor0]{Result}

\newtheorem{conj0}[theor0]{Conjecture}

\newtheorem{heur0}[theor0]{Heuristics}

\theoremstyle{definition}
\newtheorem{defin0}[theor0]{Definition}
\newenvironment{defin}
  {\pushQED{\qed}\begin{defin0}}{\popQED\end{defin0}}
\newtheorem{rems0}[theor0]{Remarks}

\newtheorem{ex0}[theor0]{Example}

\newtheorem{exs0}[theor0]{Examples}

\newtheorem{rem0}[theor0]{Remark}
\newenvironment{rem}
  {\pushQED{\qed}\begin{rem0}}{\popQED\end{rem0}}
\newtheorem{qu0}[theor0]{Question}

\newtheorem{qus0}[theor0]{Questions}

  \newtheorem{as0}[theor0]{Assumption}

\newcommand{\e}{\varepsilon}
\newcommand{\Q}{\mathbb Q}
\newcommand{\R}{\mathbb R}
\newcommand{\Z}{\mathbb Z}

\newcommand{\B}{\mathcal B}

\newcommand{\Cc}{\mathcal C}

\newcommand{\G}{\mathcal G}

\newcommand{\M}{\mathcal M}

\newcommand{\A}{\mathcal A}
\newcommand{\p}{\mathbb{P}}
\newcommand{\E}{\mathbb{E}}

\newcommand{\osc}{\partial^{\mathrm{osc}}}
\newcommand{\oscd}[1]{\mathrel{\osc_{#1}}}
\newcommand{\supess}{\operatorname{sup\,ess}}
\newcommand{\supessd}[1]{\mathrel{\mathop{\supess}\limits_{#1}}}
\newcommand{\infessd}[1]{\mathrel{\mathop{\infess}\limits_{#1}}}
\newcommand{\infess}{\operatorname{inf\,ess}}
\newcommand{\Mes}{\operatorname{Mes}}
\newcommand{\loc}{{\operatorname{loc}}}
\newcommand{\Ld}{\operatorname{L}}

\newcommand{\diam}{\operatorname{diam}}
\newcommand{\3}{\operatorname{|\!|\!|}}
\newcommand{\var}[1]{\mathrm{Var}\left[#1\right]}
\newcommand{\ent}[1]{\mathrm{Ent}\!\left[#1\right]}
\newcommand{\entm}[1]{\mathrm{Ent}[#1]}
\newcommand{\varm}[1]{\mathrm{Var}[#1]}
\newcommand{\expecm}[1]{\mathbb{E}[ #1 ]}
\newcommand{\cov}[2]{\mathrm{Cov}\left[#1;#2\right]}
\newcommand{\expec}[1]{\mathbb{E}\left[ #1 \right]}

\newcommand{\pr}[1]{\mathbb{P}\left[ #1 \right]}

\newcommand{\expeC}[2]{\mathbb{E}\left[\left. #1 \,\right\|\,#2\right]}

\newcommand{\parfct}[1]{\partial_{#1}^{\operatorname{fct}}}
\newcommand{\parG}[1]{\partial_{#1}^{\operatorname{G}}}

\newcommand{\step}[1]{\noindent \textit{Step} #1.}
\newcommand{\substep}[1]{\noindent \textit{Substep} #1.}

\numberwithin{equation}{section}

\title[Multiscale functional inequalities: Concentration properties]{Multiscale functional inequalities in probability:\\ Concentration properties}
\author[M. Duerinckx]{Mitia Duerinckx}
\author[A. Gloria]{Antoine Gloria}
\address[Mitia Duerinckx]{Laboratoire de Mathématique d'Orsay, UMR 8628, Université Paris-Sud, F-91405 Orsay, France \& Universit\'e Libre de Bruxelles, Département de Mathématique, Brussels, Belgium}
\email{mduerinc@ulb.ac.be}
\address[Antoine Gloria]{Sorbonne Universit\'e, CNRS, Universit\'e de Paris, Laboratoire Jacques-Louis Lions (LJLL), F-75005 Paris, France \& Universit\'e Libre de Bruxelles, Département de Mathématique, Brussels, Belgium}
\email{gloria@ljll.math.upmc.fr}
\begin{document}
\maketitle

\begin{abstract}
In a companion article we have introduced a notion of multiscale functional inequalities for functions $X(A)$ of an ergodic stationary random field $A$ on the ambient space $\R^d$. These inequalities are
multiscale weighted versions of standard Poincaré, covariance, and logarithmic Sobolev inequalities. They hold for all the examples of fields~$A$ arising in the modelling of heterogeneous materials in the applied sciences
whereas their standard versions are much more restrictive.
In this contribution we first investigate the link between multiscale functional inequalities and more standard decorrelation or mixing properties of random fields.
Next, we show that multiscale functional inequalities imply fine concentration properties for nonlinear functions $X(A)$. This constitutes the main stochastic ingredient to the quenched large-scale regularity theory for random elliptic operators by the second author, Neukamm, and Otto, and to the corresponding quantitative stochastic homogenization results.
\end{abstract}

\setcounter{tocdepth}{1}
\tableofcontents

\section{Introduction and main results}

Functional inequalities like Poincaré, covariance, or logarithmic Sobolev inequalities are powerful tools to prove nonlinear concentration properties and central limit theorem scalings (see e.g.~Propositions~\ref{prop:scl-averages} and \ref{prop:app-concentration} below). 
As emphasized in the introduction of the companion article~\cite{DG2}, such inequalities are also of remarkable use to study partial differential equations with random coefficients, in particular to establish quantitative results in stochastic homogenization.
However, these functional inequalities only hold under very stringent assumptions on the underlying random coefficient field and do not meet the requirements of practical models of interest to the applied sciences (see e.g.~\cite{Torquato-02}).
For this reason, we introduced in~\cite{DG2} a new family of weaker variants of these inequalities, which we refer to as \emph{multiscale functional inequalities}.
The aim of the present contribution is twofold: first we relate multiscale functional inequalities to more standard mixing conditions and decorrelation properties, 
and second we establish the concentration properties that are implied by these new inequalities
in a form to be directly used in stochastic homogenization~\cite{GNO-reg,GNO-quant,DGO2}. The present contribution can thus be viewed as an intermediate step making use and adapting classical arguments in the concentration-of-measure phenomenon literature (mainly
borrowed from~\cite{Aida-Stroock-94,Ledoux-97,Bobkov-Ledoux-97,Ledoux-01,Bakry-Gentil-Ledoux-14}) towards the applications in stochastic homogenization.
That multiscale functional inequalities yield strong concentration properties will indeed not come as a surprise to the expert and the proofs of these results constitute variations around more or less standard techniques.
For motivations and specific examples of coefficient fields that satisfy multiscale functional inequalities, we refer the reader to the companion article~\cite{DG2}.
{\color{black}Examples include Poisson random inclusions with (unbounded) random radii, Voronoi or Delaunay tessellations associated with a Poisson point process, the random parking process~\cite{Penrose-01}, and Gaussian fields with arbitrary covariance function.
More generally, all the random fields of the reference textbook \cite{Torquato-02} are covered.}
For the application of the results of the present contribution in stochastic homogenization, we refer to~\cite{GNO-reg,GNO-quant,DGO2}.
We believe that these results may be useful in other contexts for the study of partial differential equations with random coefficients.

\medskip

{\bf Notation.} 
\begin{itemize}
\item $d$ is the dimension of the ambient space $\R^d$;
\item $C$ denotes various positive constants that only depend on the dimension $d$ and possibly on other controlled quantities; we write $\lesssim$ and $\gtrsim$ for $\le$ and $\ge$ up to such multiplicative constants $C$; we use the notation $\simeq$ if both relations $\lesssim$ and $\gtrsim$ hold; we add a subscript in order to indicate the dependence of the multiplicative constants on other parameters;
\item the notation $a\gg b$ stands for $a \ge  C b$ for some large enough constant $C\simeq1$;
\item $Q:=[-1/2,1/2)^d$ denotes the unit cube centered at $0$ in dimension $d$,
and for all  $x\in\R^d$ and $r>0$ we set $Q(x):=x+Q$, $Q_r:=rQ$ and $Q_r(x):=x+rQ$;
\item we use similar notation for balls, replacing $Q$ by $B$ (the unit Euclidean ball in dimension~$d$);
\item for all subsets $D$ of $\R^d$ we denote by $\fint_D$ the average of $D$;
\item $\B(\R^k)$ denotes the Borel $\sigma$-algebra on $\R^k$;
\item $\expec{\cdot}$ denotes the expectation, $\var{\cdot}$ the variance, and $\cov{\cdot}{\cdot}$ the covariance in the underlying probability space $(\Omega,\A,\p)$,
and the notation $\expec{\cdot \| \cdot}$ stands for the conditional expectation;
\item for all $a,b\in \R$, we set $a\wedge b:=\min\{a,b\}$, $a\vee b:=\max\{a,b\}$, and $a_+:=a\vee0$.
\end{itemize}

\subsection{Multiscale functional inequalities}\label{chap:spectralgaps}

We start by recalling the notion of multiscale functional inequalities introduced in~\cite{DG2}.
Let $A:\R^d\times\Omega\to\R$ be a jointly measurable random field on $\R^d$, constructed on a probability space $(\Omega,\A,\p)$.
In what follows we assume that $A$ is stationary with respect to shifts on $\R^d$, that is, for all $z\in \R^d$ the translate $A(\cdot+z,\cdot)$ has the same finite-dimensional distributions as $A$.
A Poincaré inequality for~$A$ in the probability space is a functional inequality that allows to control the variance of any $\sigma(A)$-measurable random variable $X(A)$ in terms of its local dependence on $A$, that is, in terms of some ``derivative'' of $X(A)$ with respect to local restrictions of $A$. 
For the continuum setting that we consider here, we recall two such possible notions.
\begin{itemize}
\item The {\it oscillation} $\osc$ is formally defined by
\begin{multline}
\hspace{1cm}\oscd{A,S}X(A)~:=~\supessd{A,S}X(A)-\infessd{A,S}X(A)\\
\text{``$=$''}~ \supess\left\{X(A'):A'\in\Mes(\R^d;\R),\,A'|_{\R^d\setminus S}=A|_{\R^d\setminus S}\right\} \\
-\infess\left\{X(A'):A'\in\Mes(\R^d;\R),\,A'|_{\R^d\setminus S}=A|_{\R^d\setminus S}\right\},\label{e.def-Glauber-formal}
\end{multline}
where the essential supremum and infimum are taken with respect to the measure induced by the field $A$ on the space of measurable
real functions $\Mes(\R^d;\R)$ (endowed with the cylindrical $\sigma$-algebra). This definition \eqref{e.def-Glauber-formal} of $\osc_{A,S}X(A)$ is not measurable in general, and we rather define
\[\oscd{A,S}X(A):=\M[X\|A|_{\R^d\setminus S}]+\M[-X\|A|_{\R^d\setminus S}]\]
in terms of the conditional essential supremum $\M[\cdot\|A_{\R^d\setminus S}]$ given $\sigma(A|_{\R^d\setminus S})$, as introduced in~\cite{Barron}.
Alternatively, we may simply define $\osc_{A,S}X(A)$ as the measurable envelope of~\eqref{e.def-Glauber-formal}.
These measurable choices are equivalent for the application to stochastic homogenization, and one should not worry about these measurability issues.
\smallskip\item The (integrated) {\it functional} (or Malliavin type) {\it derivative} $\parfct{}$ is the closest generalization of the usual partial derivatives commonly used in the discrete setting.
Choose an open set $M\subset \Ld^\infty(\R^d)$ containing the realizations of the random field~$A$. Given a $\sigma(A)$-measurable random variable $X(A)$ and given an extension $\tilde X:M\to\R$ of $X$, {\color{black}its Gâteaux derivative $\frac{\partial\tilde X(A)}{\partial A}$ is defined as follows, for all compactly supported perturbation $B\in\Ld^\infty(\R^d)$, almost surely,
\[\lim_{t\to0}\frac{\tilde X(A+tB)-\tilde X(A)}t=\int_{\R^d}B(x)\,\frac{\partial\tilde X(A)}{\partial A}(x)\,dx,\]
if the limit exists with $\frac{\partial\tilde X(A)}{\partial A}\in\Ld^2(\Omega,\Ld^1_\loc(\R^d))$.
(The extension $\tilde X$ is only needed to make sure that quantities like $\tilde X(A+tB)$ make sense for small $t$, while $X$ is a priori only defined on realizations of $A$. In the sequel, we will always assume that such an extension is implicitly given; this is typically the case in applications in stochastic homogenization.)
Since we are interested in the local averages of this derivative, we rather define for all bounded Borel subset $S\subset\R^d$,
\begin{align*}
\qquad&\partial_{A,S}^{\operatorname{fct}}X(A)=\int_S\Big|\frac{\partial\tilde X(A)}{\partial A}(x)\Big|dx ~\in ~\Ld^2(\Omega).
\end{align*}}
This derivative is additive with respect to the set $S$: for all disjoint Borel subsets $S_1,S_2\subset\R^d$ we have $\partial_{A,S_1\cup S_2}^{\operatorname{fct}}X(A)=\partial_{A,S_1}^{\operatorname{fct}}X(A)+\partial_{A,S_2}^{\operatorname{fct}}X(A)$.
\end{itemize}

\medskip

Henceforth we use $\tilde\partial$ to denote either $\parfct{}$ or $\osc$. 
We are in position to recall the definition of standard functional inequalities. 
\begin{defin}\label{defi:TFI}
We say that $A$ satisfies the {\it (standard) Poincaré inequality} or {\it spectral gap \emph{($\tilde \partial$-SG)} with radius $R>0$ and constant $C>0$} if for all $\sigma(A)$-measurable random variables~$X(A)$ we have
\begin{align*}
\var{X(A)}\le \,C\,\expec{\int_{\R^d}\Big(\tilde\partial_{A,B_{R}(x)}X(A)\Big)^{2}dx};
\end{align*}
it satisfies the  {\it (standard) covariance inequality \emph{($\tilde\partial$-CI)} with radius $R>0$ and constant $C>0$} 
if for all $\sigma(A)$-measurable random variables $X(A)$ and $Y(A)$ we have
\begin{align*}
\cov{X(A)}{Y(A)}
\le C\int_{\R^d}\expec{\Big(\tilde\partial_{A,B_{R}(x)}X(A)\Big)^2}^{\frac12}\expec{\Big(\tilde\partial_{A,B_{R}(x)}Y(A)\Big)^2}^{\frac12}dx;
\end{align*}
it satisfies the {\it (standard) logarithmic Sobolev inequality \emph{($\tilde\partial$-LSI)} with radius $R>0$ and constant $C>0$} if for all $\sigma(A)$-measurable random variables $Z(A)$ we have
\[\ent{Z(A)^2}\,:=\,\expec{Z(A)^2\log\frac{Z(A)^2}{\expec{Z(A)^2}}}\,\le\,C\, \expec{\int_{\R^d}\Big(\tilde\partial_{A,B_{R}(x)}Z(A)\Big)^{2}dx}.\qedhere\]
\end{defin}

Next, we recall the definition of multiscale functional inequalities as introduced in the companion article~\cite{DG2}. These are modifications of standard functional inequalities, where derivatives with respect to the field $A$ are considered on all scales and suitably weighted. Standard inequalities are recovered for compactly supported weight $\pi$.
 
\begin{defin}

Given an integrable function $\pi:\R_+\to\R_+$, we say that $A$ satisfies the {\it multiscale Poincaré inequality \emph{($\tilde \partial$-MSG)} with weight $\pi$} if for all $\sigma(A)$-measurable random variables $X(A)$ we have
\begin{align*}
\var{X(A)}\le \,\expec{\int_0^\infty \int_{\R^d}\Big(\tilde\partial_{A,B_{\ell+1}(x)}X(A)\Big)^{2}dx\,(\ell+1)^{-d}\pi(\ell)\,d\ell}.
\end{align*}
Likewise, we define the corresponding {\it multiscale covariance inequality~\emph{($\tilde\partial$-MCI)}} and the {\it multiscale logarithmic Sobolev inequality \emph{($\tilde\partial$-MLSI)}}. 
\end{defin}

\subsection{Decay of correlations and mixing}\label{subsec:mix}

We show that the weight $\pi$ in multiscale functional inequalities is related to the decay of correlations of the random field $A$.
On the one hand, the weight is shown to control the decay of the covariance function of $A$ on large distances.
On the other hand, multiscale functional inequalities are shown to imply $\alpha$-mixing with coefficient decay related to the weight (with a caveat in the case of the functional derivative).
We start with the decay of the covariance function.
\begin{prop}\label{prop:staMSG}
Let $A$ be a jointly measurable stationary random field on $\R^d$ with $\expec{|A(0)|^2}<\infty$, and let $\Cc(x):=\cov{A(0)}{A(x)}$ denote its covariance function. 
When using the derivative $\tilde\partial=\osc$, further assume that $A$ is bounded.
\begin{enumerate}[(i)]
\item If $A$ satisfies \emph{($\tilde\partial$-SG)} and if the covariance function $\Cc$ is nonnegative, then $\Cc$ is integrable.
\item If $A$ satisfies \emph{($\tilde\partial$-MSG)} with weight $\pi$ and if the covariance function $\Cc$ is nonnegative, then $\Cc$ is integrable whenever $\int_0^\infty \ell^d \pi(\ell)d\ell<\infty$. More generally, $\Cc$ satisfies
\[\int_{\R^d}(1+|x|)^{-\alpha}\Cc(x)dx\le C_\alpha\begin{cases}\int_0^\infty (\ell+1)^{d-\alpha}\pi(\ell)d\ell,&\text{if $0\le\alpha<d$;}\\
\int_0^\infty \log^2 (2+\ell)\pi(\ell)d\ell,&\text{if $\alpha=d$;}\\
\int_0^\infty \pi(\ell)d\ell,&\text{if $\alpha>d$.}
\end{cases}\]
\item If $A$ satisfies \emph{($\tilde\partial$-CI)} with radius $R+\e$ for all $\e>0$, then the range of dependence of $A$ is bounded by $2R$ (that is, for all Borel subsets $S,T\subset\R^d$ the restrictions $A|_S$ and $A|_T$ are independent whenever $d(S,T)>2R$).
\item If $A$ satisfies \emph{($\tilde\partial$-MCI)} with weight $\pi$, then the covariance function satisfies for all $x\in\R^d$,
\[|\Cc(x)|\le C\int_{\frac12(|x|-2)\vee0}^\infty \pi(\ell)d\ell.\qedhere\]
\end{enumerate}
\end{prop}
Note in particular that {\rm($\tilde\partial$-MCI)} gives much more information than {\rm($\tilde\partial$-MSG)} on the covariance function.
As shown in the companion article~\cite[Corollary~3.1]{DG2}, Proposition~\ref{prop:staMSG} is sharp: in the Gaussian case each of the necessary conditions is (essentially) sufficient.

\medskip

We turn to ergodicity and mixing properties, and investigate the relation between multiscale Poincaré inequalities and standard mixing conditions.
Let us first recall some terminology. 
{\color{black} The random field $A$ is said to be {\it ergodic} if for all $\sigma(A)$-measurable random variables $X(A)$ there holds almost surely
$$
\lim_{R\uparrow \infty} \fint_{B_R} X(A(\cdot+x))\,dx \,=\, \expec{X(A)}.
$$
It is {\it strongly mixing} if we further have for all Borel subsets $E,E'\subset\R$}
\[\pr{X(A)\in E,\,X(A(\cdot+x))\in E'}\stackrel{|x|\uparrow\infty}\longrightarrow \pr{X(A)\in E}\,\pr{X(A)\in E'}.\]
This qualitative property can be quantified into {\it strong mixing conditions}.
A classical way to measure the dependence between two sub-$\sigma$-algebras $\G_1,\G_2\subset\A$ is the following {\it $\alpha$-mixing coefficient}, first introduced by Rosenblatt~\cite{Rosenblatt-56},
\[\alpha(\G_1,\G_2):=\sup\big\{|\p[G_1\cap G_2]-\p[G_1]\p[G_2]|\,:\,G_1\in\G_1,G_2\in\G_2\big\}.\]
Applied to the random field $A$, this leads to the following measure of mixing:
for all diameters $D\in(0,\infty]$ and distances $R>0$ we set
\begin{multline}\label{e:def-alpha-mixing}
\tilde\alpha(R,D;A):=\sup\big\{\alpha(\sigma(A|_{S_1}),\sigma(A|_{S_2}))\,:\, S_1,S_2\in\B(\R^d),\,d(S_1,S_2)\ge R,\\
\,\diam(S_1),\diam(S_2)\le D\big\}.
\end{multline}
We say that the field $A$ is {\it$\alpha$-mixing} if for all diameters $D\in(0,\infty)$ we have $\tilde \alpha(R,D;A)\stackrel{R\uparrow \infty} \longrightarrow 0$.
Note that $\alpha$-mixing is the weakest of the usual strong mixing conditions (see e.g.~\cite{Doukhan-94}), although it is in general strictly stronger than qualitative strong mixing.
\begin{prop}\label{prop:sgmix}
Let $A$ be a jointly measurable stationary random field on $\R^d$.
\begin{enumerate}[(i)]
\item If $A$ satisfies \emph{($\tilde\partial$-MSG)} with integrable weight $\pi$, then $A$ is ergodic.
\item If $A$ satisfies \emph{($\tilde\partial$-MCI)} with integrable weight $\pi$, then $A$ is strongly mixing.
\item If $A$ satisfies \emph{($\tilde\partial$-MCI)} with weight $\pi$ and derivative $\tilde\partial=\osc$, then $A$ is $\alpha$-mixing with coefficient $\tilde \alpha(R,D;A)\lesssim (1+\frac{D}{R})^d\int_{R-1}^\infty \pi(\ell)d\ell$.
\item If $A$ satisfies \emph{($\tilde\partial$-CI)} with radius $R+\e>0$ for all $\e>0$, then $\tilde\alpha(r,\infty;A)=0$ for all $r>2R$.\qedhere
\end{enumerate}
\end{prop}

\subsection{Application to spatial averages}\label{subsec:average}
Although the primary aim of functional inequalities is to address concentration properties for general \emph{nonlinear} functions of the random field $A$, we illustrate (both for the non-specialist reader and for the expert in stochastic homogenization) the use of multiscale functional inequalities on the simplest functions possible, that is, on (linear) spatial averages of the random field itself (or of a possibly nonlinear local transformation hereof).

Given a jointly measurable stationary random field  $A$,  consider a
$\sigma(A)$-measurable random variable $f(A)$
that is
approximately $1$-local with respect to the field $A$ in the following sense: for all $\ell\ge0$ we assume
\begin{align}\label{eq:as-G0}
\supessd{A}\big|f(A)-\expeC{f(A)}{A|_{B_\ell}}\big| \le Ce^{-\frac\ell C}.
\end{align}
More precisely, we use the following finer notion of approximate $1$-locality: for all $x\in\R^d$ and $\ell\ge0$ we assume
\begin{align}\label{eq:as-G}
\supessd{A}\tilde\partial_{A,B_{\ell+1}(x)}f(A) \le Ce^{-\frac1C(|x|-\ell)_+}.
\end{align}
(An important particular case is when the random variable $f(A)$ is exactly $1$-local, that is, when $f(A)$ is $\sigma(A|_{B_1})$-measurable.)
Next we set $F(x):=f(A(\cdot+x))$ for all $x\in\R^d$.
{\color{black}By definition, $F$ is a stationary random field, so that $\expec{F(x)}$ does not depend on $x\in \R^d$ and is simply denoted by $\expec{F}$.} For all $L\ge0$ we then consider the random variable
\[X_L:=X_L(A):=L^{-d}\int_{\R^d} e^{-\frac{1}L|y|}(F(y)-\expec{F})\,dy,\]
that is, the spatial average of (the nonlinear approximately local transformation $F$ of) the random field $A$ at the scale $L$. Note that the results below hold in the same form if $X_L$ is replaced by $\fint_{Q_L}(F-\expec{F})$.
We start with analyzing the scaling of the variance of $X_L$.
\begin{prop}\label{prop:scl-averages}
If $A$ satisfies  \emph{($\tilde \partial$-MSG)} with integrable weight $\pi$ and derivative $\tilde \partial=\parfct{}$ or $\osc$, and if the random variable $f(A)$ satisfies~\eqref{eq:as-G}, then we have for all $L>0$,
$$
\var{X_L} \,\lesssim \, \pi_*(L)^{-1},
$$
where we define
\begin{equation*} 
\pi_*(\ell)\,:=\,\Big(\fint_{B_\ell} \int_{|x|}^\infty \pi(s)ds dx\Big)^{-1}.\qedhere
\end{equation*}
\end{prop}
\begin{rem}
If $\pi(\ell)\simeq(\ell+1)^{-1-\beta}$ for some $\beta>0$, then
\begin{equation*}
\pi_*(\ell)\,\simeq\,
\begin{cases}
(\ell+1)^{\beta},&\text{if $\beta<d$};\\
(\ell+1)^{d}\log^{-1}(2+\ell),&\text{if $\beta=d$};\\
(\ell+1)^d,&\text{if $\beta>d$}.
\end{cases}
\end{equation*}
In particular if correlations are integrable (corresponding to the case $\beta>d$), we recover the central limit theorem scaling: $\var{X_L}\lesssim \pi_*(L)^{-1} \simeq L^{-d}$ for all $L\ge1$.
\end{rem}

Next, we study the concentration properties of the spatial average $X_L$.
The following result, which is the main stochastic ingredient of \cite{GNO-reg}, shows that the scaling crucially depends on three properties: the type of multiscale functional inequality, the type of derivative, and the decay of the weight.
The proof of items~(i) and~(ii) directly follows from the general concentration results of Propositions~\ref{prop:conc} and~\ref{prop:conc-osc} below, whereas for~(iii) we need to refine the Herbst argument and exploit the form of the spatial averages.
\begin{prop}\label{prop:app-concentration}
Assume that the random variable $f(A)$ satisfies~\eqref{eq:as-G}.
\begin{enumerate}[(i)]
\item If $A$ satisfies  \emph{($\parfct{}$-MSG)} with integrable weight $\pi$ and if $\pi_*$ is defined as in Proposition~\ref{prop:scl-averages}, then for all $\delta,L>0$,
\begin{equation}\label{e.app-conc1}
\pr{X_L\ge \delta} \, \le \, \exp\big(-\tfrac\delta C \pi_*(L)^\frac12\big).
\end{equation}
If in addition $A$ satisfies~\emph{($\parfct{}$-MLSI)} with weight $\pi$, then for all $\delta,L>0$,
\begin{equation}\label{e.app-conc2}
\pr{X_L\ge \delta}\, \le \, \exp\big(-\tfrac{\delta^2} C \pi_*(L)\big).
\end{equation}
\item If $A$ satisfies \emph{($\osc$-MSG)} with weight $\pi(\ell)\lesssim (\ell+1)^{-\beta-1}$ for some $\beta>0$, then for all $\delta,L>0$,
\begin{equation}\label{e.app-conc3-0}
\pr{X_L\ge \delta} \, \le \, C e^{-\frac1C\delta}\big(1+ \delta^{-\frac{2\beta}{d}}|\log\delta|\big)\,L^{-\beta}.
\end{equation}
\item If $A$ satisfies \emph{($\osc$-MSG)} with weight $\pi(\ell)\lesssim \exp(-\frac1C \ell^\beta)$ for some $\beta>0$, then for all $\delta>0$ and $L\ge1$,
\begin{equation}\label{e.app-conc3}
\pr{X_L\ge \delta} \, \le \, \exp\big(-\tfrac {\delta\wedge\delta^2} {C} L^{\beta\wedge\frac d2}\big).
\end{equation}
If in addition $A$ satisfies~\emph{($\osc$-MLSI)} with weight $\pi(\ell)\lesssim \exp(-\frac1C \ell^\beta)$ for some $\beta>0$, then for all $\delta>0$ and $L\ge1$,
\begin{equation}\label{e.app-conc4}
\pr{X_L\ge \delta} \, \le \, \exp\big(-\tfrac {\delta\wedge\delta^2} {C} L^{\beta\wedge d}\big).\qedhere
\end{equation}
\end{enumerate}
\qedhere
\end{prop}

\begin{rem}
If the random variable $f(A)$ is further assumed a.s.\@ bounded by a deterministic constant $C_0\ge1$, then there holds $\pr{|X_L|>C_0}=0$, hence we may replace $\delta\wedge\delta^2$ by $\frac1{C_0}\delta^2$ in~\eqref{e.app-conc3} and~\eqref{e.app-conc4}.
\end{rem}

In the case of a super-algebraic weight, it is instructive to compare these concentration results to the corresponding results implied by the $\alpha$-mixing properties of the field $A$ (see also~\cite[Appendix~A]{Armstrong-Mourrat-16}).
(Since we are mostly interested in the scaling in $L$, we do not try to optimize the $|\!\log \delta|$-dependence below.)
\begin{prop}\label{prop:app-concentration-mixing}
Given $\beta>0$, assume that the random field $A$ either is $\alpha$-mixing with coefficient $\tilde\alpha(\ell,D;A)\lesssim(1+D)^C\exp(-\frac1C\ell^\beta)$ for all $D,\ell\ge0$, or satisfies~\emph{($\osc{}$-MCI)} with weight $\pi(\ell)\lesssim\exp(-\frac1C\ell^\beta)$.
Further assume that the random variable $f(A)$ is a.s.\@ bounded by a deterministic constant and that it satisfies~\eqref{eq:as-G0}. Then for all $\delta>0$ and $L\ge1$,
\[\pr{X_L>\delta}\le C\exp\Big(-\tfrac1C{\delta^2(|\!\log \delta|+1)^{-\frac{d\beta}{d+\beta}}}\,L^{\frac{d\beta}{d+\beta}}\Big).\qedhere\]
\end{prop}

Let us briefly compare the concentration results of Propositions~\ref{prop:app-concentration}(iii) and~\ref{prop:app-concentration-mixing}. Assume that the random field $A$ satisfies a multiscale functional inequality with super-algebraic weight $\pi(\ell)\lesssim \exp(-\frac1C\ell^\beta)$ and derivative $\osc$ and that it is $\alpha$-mixing with coefficient $\tilde\alpha(\ell,D;A)\lesssim (1+D)^d\exp(-\frac1C\ell^\beta)$ (these assumptions are indeed compatible in view of Proposition~\ref{prop:sgmix}(iii)). Then, the decay in $L$ of the probability $\pr{X_L\ge\delta}$ obtained from the $\alpha$-mixing is stonger than the one obtained from ($\osc$-MSG) only for $\beta>d$, and is always weaker than the one obtained from ($\osc$-MLSI).

\subsection{Nonlinear concentration}\label{subsec:concentration}
Similarly as their standard counterparts, multiscale functional inequalities imply fine concentration properties for general random variables $X(A)$. We no longer focus on (linear) spatial averages here and consider general nonlinear measurable functions of $A$.
We start with a control on higher moments, which is a useful tool for applications (see e.g.~\cite{GNO-quant,DGO2}).
Note that the position of the weight in the right-hand side differs whether the derivative is the functional derivative or the oscillation: in particular, the control of higher moments is much weaker for the latter.
\begin{prop}\label{prop:pwsg}
Assume that $A$ satisfies \emph{($\tilde \partial$-MSG)} with integrable weight $\pi:\R_+\to\R_+$. Then there exists $C>0$ (depending only on $d,\pi$) such that for all $1\le p<\infty$ and $\sigma(A)$-measurable random variables $X(A)$,
\begin{enumerate}[(i)]
\item if $\tilde\partial=\parfct{}$, 
\begin{multline*}
\qquad \expec{\big(X(A)-\expec{X(A)}\big)^{2p}}\\
\le (Cp^{2})^p\,\expec{\bigg(\int_0^\infty\int_{\R^d}\Big(\parfct{A,B_{\ell+1}(x)}X(A)\Big)^2dx\,(\ell+1)^{-d}\pi(\ell)d\ell\bigg)^p},
\end{multline*}
where the multiplicative factor $(Cp^{2})^p$ can be upgraded to $(Cp)^p$ if the field $A$ further satisfies \emph{($\tilde \partial$-MLSI)};
\item if $\tilde\partial=\osc$,
\begin{align*}
&~\expec{\big(X(A)-\expec{X(A)}\big)^{2p}}
\\
&\hspace{1.8cm}\le (Cp^{2})^p\,\expec{\int_0^\infty\bigg(\int_{\R^d}\Big(\osc_{A,B_{2(\ell+1)}(x)}X(A)\Big)^2dx\bigg)^p{(\ell+1)^{-dp}}{\pi(\ell)}d\ell}.\qedhere
\end{align*}
\end{enumerate}
\end{prop}

From standard arguments we then deduce concentration for ``Lipschitz'' functions of the field $A$.
We start with the case of the functional derivative.
\begin{prop}\label{prop:conc}
Assume that $A$ satisfies \emph{($\parfct{}$-MSG)} with integrable weight $\pi:\R_+\to\R_+$.
We define the Lipschitz norm of a $\sigma(A)$-measurable random variable~$X(A)$ with respect to the derivative $\parfct{}$ and the weight $\pi$ as
\begin{align*}
\3X\3_{\parfct{},\pi}:=&~\supessd{A}\bigg(\int_0^\infty\int_{\R^d}\Big(\parfct{A,B_{\ell+1}(x)}X(A)\Big)^2dx\,(\ell+1)^{-d}\pi(\ell)d\ell\bigg)^{\frac12}.
\end{align*}
There exists a constant $C>0$ (depending only on $d,\pi$) such for all $\sigma(A)$-measurable random variables $X(A)$
with $\3X\3_{\parfct{},\pi}\le1$ we have exponential tail concentration in the form
\begin{gather*}
\expec{\exp\big(\tfrac1C|X(A)-\expec{X(A)}|\big)}\le2,\\
\pr{X(A)- \expec{X(A)}\ge r}\le e^{-\frac{r}{C}},\qquad \text{for all $r\ge0$}.
\end{gather*}
If in addition $A$ satisfies \emph{($\parfct{}$-MLSI)} with weight $\pi$, then for all $\sigma(A)$-measurable random variables $X(A)$ with $\3X\3_{\parfct{},\pi}\le1$ we have Gaussian tail concentration in the form
\begin{gather*}
\expec{\exp\big(\tfrac1C|X(A)-\expec{X(A)}|^2\big)}\le2,\\
\pr{X(A)- \expec{X(A)}\ge r}\le e^{-\frac{r^2\vee r}{C}},\qquad \text{for all $r\ge0$}.\qedhere
\end{gather*}
\end{prop}
We turn to the case of the oscillation, which only yields weaker concentration results due to the failure of the Leibniz rule.
\begin{prop}\label{prop:conc-osc}$  $
\begin{enumerate}[(i)]
\item Assume that $A$ satisfies \emph{($\osc$-SG)} with radius $R>0$.
Then for all $\sigma(A)$-measurable random variables $X(A)$ with
$$
\3X\3_{\osc,R}:=~\supessd{A} \int_{\R^d}\Big(\oscd{A,B_{R}(x)}X(A)\Big)^2dx\,\le \,1,
$$
we have exponential tail concentration in the form
\begin{gather*}
\expec{\exp\big(\tfrac1C|X(A)-\expec{X(A)}|\big)}\le2,\\
\pr{X(A)- \expec{X(A)}\ge r}\le e^{-\frac{r}{C}},\qquad \text{for all $r\ge0$}.
\end{gather*}
If in addition $A$ satisfies \emph{($\osc$-LSI)} with radius $R>0$ and if the random variable $X(A)$ further satisfies
\[L:=\sup_{x}\supessd{A}\osc_{A,B_{R}(x)}X(A)<\infty,\]
we have Poisson tail concentration in the form
\begin{gather*}
\expec{\exp\big(\tfrac1C\psi_L(|X(A)-\expec{X(A)}|)\big)}\le 2,\qquad \psi_L(u):=\tfrac{u}{L}\log\big(1+\tfrac{L}Cu\big),\\
\pr{X(A)- \expec{X(A)}\ge r}\le e^{-\frac 1{C}\psi_L(r)},\qquad \text{for all $r\ge0$}.
\end{gather*}
\item Assume that $A$ satisfies~\emph{($\osc$-MSG)} with integrable weight $\pi:\R_+\to \R_+$.
Let $X(A)$ be a $\sigma(A)$-measurable random variable and assume that for some $\kappa>0$, $p_0,\alpha\ge 0$ we have for all $p\ge p_0$,
\begin{align}\label{e.conc-osc}
\expec{\int_0^\infty\bigg(\int_{\R^d}\Big(\oscd{A,B_{\ell+1}(x)}X(A)\Big)^{2}dx\bigg)^p\,(\ell+1)^{-dp}\pi(\ell)d\ell} \,\le \, p^{\alpha p}\kappa.
\end{align}
Then there is a constant $C>0$ (depending only on $d,\pi,p_0,\alpha$ but not on $\kappa$) such that we have concentration in the form 
\begin{gather*}
\qquad\expec{\psi_{p_0,\alpha}\big(\tfrac1C|X(A)-\expec{X(A)}|\big)}\,\le\,C\kappa,\qquad \psi_{p_0,\alpha}(u):=(1\wedge r^{2p_0})\exp(r^\frac{2}{2+\alpha}),\\
\quad\pr{\,|X(A)- \expec{X(A)}|\ge r}\le C\kappa \big(\psi_{p_0,\alpha}(\tfrac rC)\big)^{-1},\qquad \text{for all $r\ge0$}.\qedhere
\end{gather*}
\end{enumerate}
\end{prop}

\medskip
\section{Decay of correlations and mixing} \label{chap:spectralgaps}
 This section is devoted to the proof of the results stated in Section~\ref{subsec:mix}, that are, Propositions~\ref{prop:staMSG} and~\ref{prop:sgmix}.

\subsection{Proof of Proposition~\ref{prop:staMSG}: Decay of correlations}

We split the proof into four steps.

\medskip

\step1 Proof of~(i).

\nopagebreak
Let the field $A$ satisfy ($\tilde\partial$-SG) with radius $R$. For any $L\ge1$, the standard Poincaré inequality applied to the $\sigma(A)$-measurable random variable $X(A)=\int_{B_L}A$ (which is well-defined by measurability and moment bounds on $A$) yields
\[\var{\int_{B_L}A}\le C\expec{\int_{\R^d}\Big(\tilde\partial_{A,B_R(x)}\int_{B_L}A\Big)^2dx}.\]
For each choice of the derivative $\tilde\partial$ (further assuming that $A$ is bounded in the case $\tilde\partial=\osc$), we have
\[\expec{\Big(\tilde\partial_{A,B_R(x)}\int_{B_L}A\Big)^2}\le C|B_R(x)\cap B_L|^2\le C_R\mathds1_{|x|\le R+L}.\]
Hence, for $L\ge1$,
\[ \int_{B_L}\int_{B_L} \cov{A(x)}{A(y)}dxdy=\var{\int_{B_L}A}\le C_R|B_{R+L}|\le C_R|B_L|.\]
Therefore, if $\Cc$ is nonnegative, we deduce
\[\int_{B_{L}}\Cc\lesssim \int_{B_L}\fint_{B_L}\Cc(x-y)dydx=\int_{B_L}\fint_{B_L}\cov{A(x)}{A(y)}dydx \le C_R.\]
Letting $L\uparrow\infty$, we conclude that $\Cc$ is integrable.

\medskip

\step2 Proof of~(ii).

\nopagebreak
Let the field $A$ satisfy ($\tilde\partial$-MSG) with weight $\pi$, and assume that $\Cc$ is nonnegative.
Repeating the argument of Step~1, we deduce for all $L\ge1$,
\begin{align*}
L^d\int_{B_L}\Cc(x)dx~&\lesssim~ \expec{\Big(\int_{B_L}(A(x)-\expec{A})dx\Big)^2}
\\
&\le~  \int_0^\infty \int_{\R^d}|B_{\ell+1}(x)\cap B_L|^2dx\,(\ell+1)^{-d}\pi(\ell)d\ell\\
&\lesssim~ \int_0^L L^d(\ell+1)^{d}\pi(\ell)d\ell+\int_L^\infty  L^{2d}\pi(\ell)d\ell \\
&\lesssim~ L^d\int_0^\infty (\ell+1)^{d}\pi(\ell)d\ell,
\end{align*}
which shows that $\Cc$ is integrable if $\int_0^\infty (\ell+1)^{d}\pi(\ell)d\ell<\infty$.

Let now $\alpha>0$ be fixed, and let $\gamma:=\frac12(d+\alpha)$. Assume that $\alpha\ne d$ (the case $\alpha=d$ can be treated similarly
and yields the logarithmic correction). For all $L\ge1$, the multiscale Poincaré inequality applied to the $\sigma(A)$-measurable random variable $X(A)=\int_{B_L} (1+|y|)^{-\gamma}A(y)dy$ yields
\begin{eqnarray*}
\lefteqn{\var{\int_{B_L} (1+|y|)^{-\gamma}A(y)dy}}\\
&\le& \expec{\int_0^\infty\int_{\R^d}\Big(\tilde\partial_{A,B_{\ell+1}(x)}\int_{B_L} (1+|y|)^{-\gamma}A(y)dy\Big)^2dx\,(\ell+1)^{-d}\pi(\ell)d\ell}\\
&\le &C\int_0^\infty\int_{\R^d}\Big(\int_{B_L\cap B_{\ell+1}(x)} (1+|y|)^{-\gamma}dy\Big)^2dx\,(\ell+1)^{-d}\pi(\ell)d\ell.
\end{eqnarray*}
Hence,
\begin{multline*}
\int_{B_{2L}}\bigg(\int_{B_L(-x)}(1+|x+y|)^{-\gamma}(1+|y|)^{-\gamma}dy\bigg)\Cc(x)dx\\
=\var{\int_{B_L} (1+|y|)^{-\gamma}A(y)dy}\le C_\alpha\int_0^\infty (\ell+1)^{(d-\alpha)\vee0}\pi(\ell)d\ell,
\end{multline*}
which yields the claim by passing to the limit $L\uparrow\infty$.

\medskip

\step3 Proof of~(iii).

Let the field $A$ satisfy ($\tilde\partial$-CI) with radius $R+\e$ for any $\e>0$.
Given two Borel subsets $S,T\subset\R^d$ with $d(S,T)>2R$, choosing $\e:=\frac13(d(S,T)-2R)$, and noting that the sets $S+B_{R+\e}$ and $T+B_{R+\e}$ are disjoint, the covariance inequality~($\tilde\partial$-CI) with radius $R+\e$ implies for any $G\in\sigma(A|_S)$ and $H\in\sigma(A|_T)$,
\begin{multline*}
|\cov{\mathds1_G}{\mathds1_H}|
\\
\le C_\e\int_{(S+B_{R+\e})\cap(T+B_{R+\e})} \expec{\Big(\tilde\partial_{A,B_{R+\e}(x)}\mathds1_G)\Big)^2}^{\frac12}\expec{\Big(\tilde\partial_{A,B_{R+\e}(x)}\mathds1_H\Big)^2}^{\frac12}dx~=0,
\end{multline*}
{\color{black}hence $\pr{G\cap H}=\pr{G}\pr{H}$, showing that the $\sigma$-algebras $\sigma(A|_S)$ and $\sigma(A|_T)$ are independent.}

\medskip

\step4 Proof of~(iv).

\nopagebreak
Let the field $A$ satisfy ($\tilde\partial$-MCI) with weight $\pi$. For all $x\in\R^d$ and all $\e>0$, the covariance inequality applied to the $\sigma(A)$-measurable random variables $\fint_{B_\e(x)} A$ and $\fint_{B_\e} A$
yields 
\begin{eqnarray*}
\lefteqn{\bigg|\fint_{B_\e(x)}\fint_{B_\e}\Cc(y-z)dydz\bigg|=\bigg|\cov{\fint_{B_\e(x)}A}{\fint_{B_\e}A}\bigg|}\\
&\le& \int_{0}^\infty\int_{\R^d}\expec{\Big(\tilde\partial_{A,B_{\ell+1}(y)}\fint_{B_\e(x)}A\Big)^2}^\frac12\expec{\Big(\tilde\partial_{A,B_{\ell+1}(y)}\fint_{B_\e}A\Big)^2}^\frac12dy\,(\ell+1)^{-d}\pi(\ell)d\ell\\
&\le& \int_{0}^\infty\int_{\R^d} \e^{-d}|B_\e(x)\cap B_{\ell+1}(y)|\e^{-d}|B_\e\cap B_{\ell+1}(y)|dy\,(\ell+1)^{-d}\pi(\ell)d\ell.
\end{eqnarray*}
Letting $\e\downarrow0$ and using the continuity of the function $\Cc$ (as a consequence of the stochastic continuity of the field $A$, which follows from its joint measurability), we 
deduce for all $x\in\R^d$,
\[|\Cc(x)|\le C\int_{0}^\infty |B_{\ell+1}(x)\cap B_{\ell+1}|\,(\ell+1)^{-d}\pi(\ell)d\ell\,\le\, C\int_{\frac12(|x|-2)\vee 0}^\infty  \pi(\ell)d\ell,\]
and the claim follows.\qed

\subsection{Proof of Proposition~\ref{prop:sgmix}: Mixing}\label{chap:erg-mix}

Item~(iv) follows from Proposition~\ref{prop:staMSG}. We split the rest of the proof into three steps.

\medskip

\step1 Proof of (i).

Let the field $A$ satisfy ($\tilde\partial$-MSG) with weight $\pi$. To prove ergodicity, it suffices to show that for all integrable $\sigma(A)$-measurable random variables $X(A)$ we have
\[\lim_{L\uparrow\infty}\expec{\Big|\fint_{B_L}X(A(x+\cdot))dx-\expec{X(A)}\Big|}=0.\]
By an approximation argument in $\Ld^2(\Omega)$, we may assume that $X(A)$ is bounded and is $\sigma(A|_{B_R})$-measurable for some $R>0$. The Poincaré inequality ($\tilde\partial$-MSG) applied
to the $\sigma(A)$-measurable random variable $\fint_{B_L}X(A(\cdot+x))dx$ yields
\begin{multline*}
S_L:=\expec{\Big|\fint_{B_L}X(A(x+\cdot))dx-\expec{X(A)}\Big|}^2\le\var{\fint_{B_L}X(A(x+\cdot))dx}\\
\le \expec{\int_0^\infty\int_{\R^d}\Big(\fint_{B_L}\tilde\partial_{A,B_{\ell+1}(y)}X(A(x+\cdot))dx\Big)^2dy\,(\ell+1)^{-d}\pi(\ell)d\ell},
\end{multline*}
and therefore
\begin{multline*}
S_L\,\le\,  \E\bigg[\int_0^\infty\int_{\R^d}\fint_{B_L}\fint_{B_L}\tilde\partial_{A,B_{\ell+1}(y)}X(A(x+\cdot))\,\tilde\partial_{A,B_{\ell+1}(y)}X(A(x'+\cdot))dxdx'dy\\
\times (\ell+1)^{-d}\pi(\ell)d\ell\bigg].
\end{multline*}
By assumption, $\tilde\partial_{A,B_{\ell+1}(y)}X(A(x+\cdot))=0$ whenever $B_R(x)\cap B_{\ell+1}(y)=\varnothing$, i.e.\@ whenever $|x-y|>R+\ell+1$. For the choices $\tilde\partial=\osc$ and $\parG{}$, we also have $\tilde\partial_{A,B_{\ell+1}(y)}X(A(x+\cdot))\le 2\|X\|_{\Ld^\infty}$, so that the above yields
\begin{align*}
S_L&\le4 \|X\|_{\Ld^\infty}^2\int_0^\infty\int_{\R^d}\fint_{B_L}\fint_{B_L}\mathds1_{|x-y|\le R+\ell+1}\mathds1_{|x'-y|\le R+\ell+1}dxdx'dy\,(\ell+1)^{-d}\pi(\ell)d\ell\\
&= 4\|X\|_{\Ld^\infty}^2L^{-2d}\int_0^\infty \Big(\int_{B_L}\int_{B_{R+\ell+1}(x)}|B_L\cap B_{R+\ell+1}(y)|dydx\Big)(\ell+1)^{-d}\pi(\ell)d\ell\\
&\le4 \|X\|_{\Ld^\infty}^2\int_0^\infty (R+\ell+1)^d  \Big(\frac{R+\ell}{L}\wedge1\Big)^d(\ell+1)^{-d} \pi(\ell)d\ell,
\end{align*}
where the RHS obviously goes to $0$ as $L\uparrow\infty$ whenever $\int_0^\infty \pi(\ell)d\ell<\infty$. This proves ergodicity for the choice $\tilde\partial=\osc$.

It remains to treat the case $\tilde\partial=\parfct{}$. An additional approximation argument is then needed in order to restrict attention to those random variables $X(A)$ such that the derivative $\tilde\partial_{A,B_{\ell+1}(x)}X(A)$ is pointwise bounded.
The stochastic continuity of the field $A$ (which follows from its joint measurability) ensures that the $\sigma(A|_{B_R})$-measurable random variable $X(A)$ is actually $\sigma(A|_{\Q^d\cap B_R})$-measurable. A standard approximation argument then allows to construct a sequence $(x_n)_n\subset B_R$ and a sequence $(X_n(A))_n$ of random variables such that $X_n(A)$ is $\sigma((A(x_k))_{k=1}^n)$-measurable and converges to $X(A)$ in $\Ld^2(\Omega)$. By definition, we may write $X_n(A)=f_n(A(x_1),\ldots,A(x_n))$ for some Borel function $f_n:(\R^k)^n\to\R$. Another standard approximation argument now allows to replace the Borel maps $f_n$'s by smooth functions. We end up with a sequence that approximates $X(A)$ in $\Ld^2(\Omega)$, and such that the elements have pointwise bounded $\tilde\partial$-derivative.
For these approximations, the conclusion follows as before.

\medskip

\step2 Proof of (ii).

Let the field $A$ satisfy ($\tilde\partial$-MCI) with weight $\pi$. To prove strong mixing, it suffices to show that for all bounded $\sigma(A)$-measurable random variables $X(A)$ and $Y(A)$ we have $\cov{X(A)}{Y(A(x+\cdot))}\to 0$ as $|x|\to\infty$ (since the desired property then follows by choosing the random variables $X(A),Y(A)$ to be any pair of indicator functions). 
Again, a standard approximation argument allows one to consider bounded $\sigma(A|_{B_R})$-measurable random variables $X(A),Y(A)$ for some $R>0$. Given $x\in\R^d$, apply the covariance inequality ($\tilde\partial$-MCI) to $X(A)$ and $Y(A(\cdot+x))$ to obtain
\begin{multline*}
\big|\cov{X(A)}{Y(A(x+\cdot))}\big|\\
\le \int_0^\infty\int_{\R^d} \expec{\Big(\tilde\partial_{A,B_{\ell+1}(y)}X(A)\Big)^2}^{\frac12}\expec{\Big(\tilde\partial_{A,B_{\ell+1}(y)}Y(A(x+\cdot))\Big)^2}^{\frac12}dy\,(\ell+1)^{-d}\pi(\ell)d\ell.
\end{multline*}
By assumption, $\tilde\partial_{A,B_{\ell+1}(y)}X(A)=0$ whenever $B_R\cap B_{\ell+1}(y)=\varnothing$, i.e.\@ whenever $|y|>R+\ell+1$. For $\tilde\partial=\osc$, we have in addition $\tilde\partial_{A,B_{\ell+1}(y)}X(A)\le 2\|X\|_{\Ld^\infty}$, so that the above directly yields
\begin{eqnarray*}
\lefteqn{\big|\cov{X(A)}{Y(A(x+\cdot))}\big|}
\\
&\le& 4\|X\|_{\Ld^\infty}\|Y\|_{\Ld^\infty}\int_{0}^\infty\int_{\R^d} \mathds1_{|y|\le R+\ell+1}\mathds1_{|x-y|\le R+\ell+1}\,dy\,(\ell+1)^{-d}\pi(\ell)d\ell\\
&\lesssim & \|X\|_{\Ld^\infty}\|Y\|_{\Ld^\infty}\int_0^\infty (R+\ell+1)^d (\ell+1)^{-d}\pi(\ell)d\ell
\end{eqnarray*}
where the RHS goes to $0$ as $|x|\to\infty$ whenever $\int_0^\infty \pi(\ell)d\ell<\infty$. This proves strong mixing for the choice $\tilde\partial=\osc$. In the case $\tilde\partial=\parfct{}$, an additional approximation argument is needed as in Step~1 in order to restrict to random variables $X(A)$ such that $\tilde\partial_{A,B_{\ell+1}(y)}X(A)$ is pointwise bounded.

\medskip

\step3 Proof of (iii).

Let the field $A$ satisfy ($\tilde\partial$-MCI) with weight $\pi$, and with derivative $\osc$. Given Borel subsets $S,T\subset\R^d$ with diameter $\le D$ and with $d(S,T)\ge 2R$, the covariance inequality ($\tilde\partial$-MCI) for this choice of derivatives yields for all bounded random variables $X(A)$ and $Y(A)$, respectively $\sigma(A|_{S})$-measurable and $\sigma(A|_{T})$-measurable,
\begin{eqnarray*}
\lefteqn{\big|\cov{X(A)}{Y(A)}\big|}
\\
&\le& \int_0^\infty \int_{\R^d}\expec{\Big(\osc_{A,B_{\ell+1}(x)}X(A)\Big)^2}^\frac12\expec{\Big(\osc_{A,B_{\ell+1}(x)}Y(A)\Big)^2}^\frac12dx\,(\ell+1)^{-d}\pi(\ell)d\ell\\
&\le& 4\|X(A)\|_{\Ld^\infty} \|Y(A)\|_{\Ld^\infty}\int_0^\infty \big|(S+B_{\ell+1})\cap(T+B_{\ell+1})\big|\,(\ell+1)^{-d}\pi(\ell)d\ell\\
&\lesssim& \|X(A)\|_{\Ld^\infty} \|Y(A)\|_{\Ld^\infty}\int_{R-1}^\infty (\ell+D+1)^d(\ell+1)^{-d}\pi(\ell)d\ell
\\
&\le &  \|X(A)\|_{\Ld^\infty} \|Y(A)\|_{\Ld^\infty}\Big(1+\frac{D}{R}\Big)^d\int_{R-1}^\infty \pi(\ell)d\ell,
\end{eqnarray*}
from which the claim follows by choosing for $X(A),Y(A)$ any indicator functions.
\qed


\section{Nonlinear concentration}\label{sec:concentration}

 This section is devoted to the proof of the results stated in Section~\ref{subsec:concentration}, that are, Proposition~\ref{prop:pwsg} and~\ref{prop:conc-osc}. The proof of Proposition~\ref{prop:conc} is obvious based on Proposition~\ref{prop:pwsg} and is omitted.

\subsection{Proof of Proposition~\ref{prop:pwsg}: Higher moments}

Let $X(A)$ be $\sigma(A)$-measurable, and let $\tilde \partial$ denote $\osc{}$ or $\parfct{}$. We may assume without loss of generality that $\expec{X(A)}=0$.
We split the proof into two steps.

\medskip

\step1 Proof of (i) and (ii) for {\rm($\tilde \partial$-MSG)}.

Applying the Poincaré inequality ($\tilde\partial$-MSG) to the $\sigma(A)$-measurable random variable $|X(A)|^p$ yields
\begin{multline}\label{eq:sg-ppower}
\expec{X(A)^{2p}}\le \expec{|X(A)|^p}^2\\
+\expec{\int_0^\infty\int_{\R^d}\Big(\tilde\partial_{A,B_{\ell+1}(x)}\big(|X(A)|^p\big)\Big)^2dx\,(\ell+1)^{-d}\pi(\ell)d\ell}.
\end{multline}
For $p>2$, H\"older's and Young's inequalities with exponents $(\frac{2(p-1)}{p-2},\frac{2(p-1)}{p})$ and $(\frac{p-1}{p-2},p-1)$, respectively, imply for all $\delta>0$,
\begin{align*}
\expec{|X(A)|^p}^2=\expec{|X(A)|^{p\frac{p-2}{p-1}}|X(A)|^{\frac{p}{p-1}}}^2&\le\expec{X(A)^{2p}}^{\frac{p-2}{p-1}}\expec{X(A)^2}^{\frac p{p-1}}\\
&\le \frac{p-2}{p-1}\delta\, \expec{X(A)^{2p}}+\frac1{p-1}\delta^{2-p}\,\expec{X(A)^2}^p.
\end{align*}
while for $p\le2$ Jensen's inequality simply yields $\expec{|X(A)|^p}^2\le\expec{X(A)^2}^p$. 
Injecting these estimates into~\eqref{eq:sg-ppower} for some $\delta\gtrsim1$ small enough, we conclude for all $1\le p<\infty$,
\begin{multline*}
\expec{X(A)^{2p}}\le p^{-1}C^p\,\expec{X(A)^2}^p\\
+C\,\expec{\int_0^\infty\int_{\R^d}\Big(\tilde\partial_{A,B_{\ell+1}(x)}\big(|X(A)|^p\big)\Big)^2dx\,(\ell+1)^{-d}\pi(\ell)d\ell}.
\end{multline*}
Since $\expec{X(A)^2}=\var{X(A)}$ follows from the centering assumption, the first RHS term is estimated by the Poincaré inequality ($\tilde\partial$-MSG). Further using Jensen's inequality, this leads to
\begin{multline}\label{eq:sg-ppower1}
\expec{X(A)^{2p}}\le p^{-1}C^p\,\expec{\bigg(\int_0^\infty\int_{\R^d}\Big(\tilde\partial_{A,B_{\ell+1}(x)}X(A)\Big)^2dx\,(\ell+1)^{-d}\pi(\ell)d\ell\bigg)^p}\\+C\,\expec{\int_0^\infty\int_{\R^d}\Big(\tilde\partial_{A,B_{\ell+1}(x)}\big(|X(A)|^p\big)\Big)^2dx\,(\ell+1)^{-d}\pi(\ell)d\ell}.
\end{multline}
We split the rest of this step into two further substeps, and treat separately $\parfct{}$ and~$\osc$.

\medskip

\substep{1.1} Proof of~(i) for $\tilde\partial=\parfct{}$.

\nopagebreak
By the Leibniz rule, $\partial_{A,S}^{\operatorname{fct}}(|X(A)|^p)=p|X(A)|^{p-1}\partial_{A,S}^{\operatorname{fct}} X(A)$, so
that H\"older's inequality with exponents $(\frac{p}{p-1},p)$ yields
\begin{eqnarray}
\lefteqn{\expec{\int_0^\infty\int_{\R^d}\Big(\partial^{\operatorname{fct}}_{A,B_{\ell+1}(x)}\big(|X(A)|^p\big)\Big)^2dx\,(\ell+1)^{-d}\pi(\ell)d\ell}}\nonumber
\\
&\le& p^2\,\expec{X(A)^{2(p-1)}\int_0^\infty\int_{\R^d}\Big(\partial^{\operatorname{fct}}_{A,B_{\ell+1}(x)}X(A)\Big)^{2}dx\,(\ell+1)^{-d}\pi(\ell)d\ell}\nonumber
\\
&\le& p^2\,\expec{X(A)^{2p}}^{1-\frac1p}\expec{\bigg(\int_0^\infty\int_{\R^d}\Big(\partial^{\operatorname{fct}}_{A,B_{\ell+1}(x)}X(A)\Big)^{2}dx\,(\ell+1)^{-d}\pi(\ell)d\ell\bigg)^p}^\frac1p.
\qquad \label{e.add-1}
\end{eqnarray}
Combined with~\eqref{eq:sg-ppower1} and Young's inequality with exponents $(\frac{p}{p-1},p)$ to absorb the factor $\expec{X(A)^{2p}}$ into the LHS, the conclusion of item~(i) follows with the prefactor $(Cp^2)^p$.

\medskip

\substep{1.2} Proof of~(ii) for $\tilde\partial=\osc$.

The inequality $||a|^p-|b|^p|\le p|a-b|(|a|^{p-1}+|b|^{p-1})$ for all $a,b\in\R$ implies
\begin{align}
\oscd{A,S}|X(A)|^p &\le 2p\Big(\sup_{A,S}|X(A)|^{p-1}\Big)\oscd{A,S}X(A) \nonumber
\\
&\le2p\bigg(|X(A)|+\oscd{A,S}X(A)\bigg)^{p-1}\oscd{A,S}X(A).
\label{eq:boundoscLp-pre}
\end{align}
We then make use of the following inequality that holds for some constant $C\simeq 1$ large enough (independent of $p$): for all $a,b\ge0$,
$(a+b)^{p-1}\le 2a^{p-1}+(Cp)^{p}b^{p-1}$. This allows one to rewrite \eqref{eq:boundoscLp-pre} in
the form
\begin{align}\label{eq:boundoscLp}
\oscd{A,S}|X(A)|^p&\le 4p|X(A)|^{p-1}\oscd{A,S}X(A)+(Cp)^p(\oscd{A,S}X(A))^p.
\end{align}
Arguing as in Substep~1.1, we obtain by H\"older's inequality,
\begin{multline*}
\expec{\int_0^\infty\int_{\R^d}\Big(\oscd{A,B_{\ell+1}(x)}|X(A)|^p\Big)^2dx\,(\ell+1)^{-d}\pi(\ell)d\ell}
\\
\le Cp^2\,\expec{X(A)^{2p}}^{1-\frac1p}\expec{\bigg(\int_0^\infty\int_{\R^d}\Big(\oscd{A,B_{\ell+1}(x)}X(A)\Big)^{2}dx\,(\ell+1)^{-d}\pi(\ell)d\ell\bigg)^p}^\frac1p\\
+(Cp^2)^{p}\,\expec{\int_0^\infty\int_{\R^d}\Big(\oscd{A,B_{\ell+1}(x)}X(A)\Big)^{2p}dx\,(\ell+1)^{-d}\pi(\ell)d\ell}.
\end{multline*}
Combined with~\eqref{eq:sg-ppower1} and Young's inequality to absorb the factor $\expecm{X(A)^{2p}}$ into the LHS, 
this yields
\begin{multline*}
\expec{X(A)^{2p}}\le (Cp^{2})^p\,\expec{\bigg(\int_0^\infty\int_{\R^d}\Big(\oscd{A,B_{\ell+1}(x)}X(A)\Big)^2dx\,(\ell+1)^{-d}\pi(\ell)d\ell\bigg)^p}\\
+(Cp^2)^{p}\,\expec{\int_0^\infty\int_{\R^d}\Big(\oscd{A,B_{\ell+1}(x)}X(A)\Big)^{2p}dx\,(\ell+1)^{-d}\pi(\ell)d\ell}.
\end{multline*}
It remains to reformulate the second RHS term. By the discrete $\ell^1-\ell^p$ inequality, we have
\begin{align}
\int_{\R^d}\Big(\oscd{A,B_{\ell+1}(x)}X(A)\Big)^{2p}dx&\le\sum_{z\in \frac{\ell+1}{\sqrt d}\Z^d}\int_{z+\frac{\ell+1}{\sqrt d}Q}\Big(\oscd{A,B_{\ell+1}(x)}X(A)\Big)^{2p}dx\nonumber\\
&\le\Big(\frac{\ell+1}{\sqrt d}\Big)^d\sum_{z\in\frac{\ell+1}{\sqrt d}\Z^d}\Big(\oscd{A,B_{\frac32(\ell+1)}(z)}X(A)\Big)^{2p}\nonumber\\
&\le\Big(\frac{\ell+1}{\sqrt d}\Big)^d\bigg(\sum_{z\in\frac{\ell+1}{\sqrt d}\Z^d}\Big(\oscd{A,B_{\frac32(\ell+1)}(z)}X(A)\Big)^{2}\bigg)^p\nonumber\\
&\le\Big(\frac{\ell+1}{\sqrt d}\Big)^d\bigg(\sum_{z\in\frac{\ell+1}{\sqrt d}\Z^d}\fint_{z+\frac{\ell+1}{\sqrt d}Q}\Big(\oscd{A,B_{2(\ell+1)}(x)}X(A)\Big)^{2}dx\bigg)^p\nonumber\\
&\le\Big(\frac{\sqrt d}{\ell+1}\Big)^{d(p-1)}\bigg(\int_{\R^d}\Big(\oscd{A,B_{2(\ell+1)}(x)}X(A)\Big)^{2}dx\bigg)^p.\label{eq:estl1lp}
\end{align}
Combined with the above, this yields
\begin{multline*}
\expec{X(A)^{2p}}\le (Cp^{2})^p\,\expec{\bigg(\int_0^\infty\int_{\R^d}\Big(\oscd{A,B_{\ell+1}(x)}X(A)\Big)^2dx\,(\ell+1)^{-d}\pi(\ell)d\ell\bigg)^p}\\
+(Cp^2)^{p}\,\expec{\int_0^\infty\bigg(\int_{\R^d}\Big(\oscd{A,B_{2(\ell+1)}(x)}X(A)\Big)^{2}dx\bigg)^p\,(\ell+1)^{-dp}\pi(\ell)d\ell}.
\end{multline*}
Since $\int_0^\infty \pi(\ell)d\ell<\infty$, the first RHS term can be absorbed into the second RHS term.
Indeed, the triangle inequality and the H\"older inequality with exponents $(p,\frac p{p-1})$ combine to
\begin{eqnarray*}
\lefteqn{\expec{\bigg(\int_0^\infty\int_{\R^d}\Big(\oscd{A,B_{\ell+1}(x)}X(A)\Big)^2dx\,(\ell+1)^{-d}\pi(\ell)d\ell\bigg)^p}}\\
&\le&\bigg(\int_0^\infty\expec{\bigg(\int_{\R^d}\Big(\oscd{A,B_{\ell+1}(x)}X(A)\Big)^2dx\bigg)^p}^{\frac1p}(\ell+1)^{-d}\pi(\ell)d\ell\bigg)^p\\
&=&\bigg(\int_0^\infty\expec{\bigg(\int_{\R^d}\Big(\oscd{A,B_{\ell+1}(x)}X(A)\Big)^2dx\bigg)^p(\ell+1)^{-dp}\pi(\ell)}^{\frac1p} \pi(\ell)^{1-\frac1p}d\ell\bigg)^p\\
&\le &\bigg(\int_0^\infty\pi(\ell)d\ell\bigg)^{p-1}\expec{\int_0^\infty\bigg(\int_{\R^d}\Big(\oscd{A,B_{\ell+1}(x)}X(A)\Big)^2dx\bigg)^p(\ell+1)^{-dp}\pi(\ell)d\ell},
\end{eqnarray*}
and the conclusion of item~(ii) follows.

\medskip

\step2 Improvement of (i) for  {\rm($\parfct{}$-MLSI)}.

\nopagebreak
In this step, we argue that the prefactor $(Cp^2)^p$ in item~(i) can be upgraded to $(Cp)^p$ if the field $A$ satisfies the corresponding logarithmic Sobolev inequality {\rm($\parfct{}$-MLSI)}.
Starting point is the following observation (see~\cite[Theorem~3.4]{Aida-Stroock-94} and~\cite[Proposition~5.4.2]{Bakry-Gentil-Ledoux-14}): if the random variable $X(A)$ satisfies $\ent{X(A)^{2p}}<\infty$, then we have
\begin{align}\label{eq:der-Lp}
\expec{X(A)^{2p}}^\frac1p-\expec{X(A)^{2}}=\int_1^p\frac1{q^2}\expec{X(A)^{2q}}^{\frac1q-1}\ent{X(A)^{2q}}dq.
\end{align}
It remains to estimate the entropy $\entm{X(A)^{2q}}$ for all $1\le q\le p$. Applied to the $\sigma(A)$-measurable random variable $|X(A)|^q$,
($\parfct{}$-MLSI) yields
\begin{align*}
\ent{X(A)^{2q}}\le \,\expec{\int_0^\infty\int_{\R^d}\Big(\parfct{A,B_{\ell+1}(x)}|X(A)|^q\Big)^2dx\,(\ell+1)^{-d}\pi(\ell)d\ell}.
\end{align*}
The argument of Substep~1.1, cf.~\eqref{e.add-1}, applied to the above RHS yields
\begin{multline*}
\ent{X(A)^{2q}}\\
\le Cq^2\,\expec{X(A)^{2q}}^{1-\frac1q}\expec{\bigg(\int_0^\infty\int_{\R^d}\Big(\parfct{A,B_{\ell+1}(x)}X(A)\Big)^{2}dx\,(\ell+1)^{-d}\pi(\ell)d\ell\bigg)^q}^\frac1q.
\end{multline*}
Inserting this into~\eqref{eq:der-Lp}, we obtain
\begin{multline*}
\expec{X(A)^{2p}}^\frac1p\le\expec{X(A)^{2}}
\\
+C\int_1^p\expec{\bigg(\int_0^\infty\int_{\R^d}\Big(\parfct{A,B_{\ell+1}(x)}X(A)\Big)^{2}dx\,(\ell+1)^{-d}\pi(\ell)d\ell\bigg)^q}^\frac1qdq.
\end{multline*}
We then appeal to the Poincaré inequality~($\parfct{}$-MSG) (which follows from ($\parfct{}$-MLSI)) to estimate the first RHS term, and use Jensen's inequality on the second RHS to obtain
\begin{align*}
\expec{X(A)^{2p}}^\frac1p\le Cp\,\expec{\bigg(\int_0^\infty\int_{\R^d}\Big(\parfct{A,B_{\ell+1}(x)}X(A)\Big)^{2}dx\,(\ell+1)^{-d}\pi(\ell)d\ell\bigg)^p}^\frac1p.
\end{align*}
This upgrades the prefactor in item~(i) to $(Cp)^p$, as claimed.\qed

\subsection{Proof of Proposition~\ref{prop:conc-osc}: Nonlinear concentration}

We split the proof into two steps, and prove (i) and (ii) separately.

\medskip
\step1 Proof of~(i).

The exponential concentration result in~(i) follows from Proposition~\ref{prop:pwsg}(ii) (with compactly supported weight $\pi$) by standard 
arguments {\color{black}(e.g.~\cite[Proposition~2.13]{Ledoux-97})}.
Let us now turn to the Poisson concentration result, and use a variation of the Herbst argument.
Let $A$ satisfy ($\osc$-LSI) and let $X(A)$ satisfy $L:=\sup_{x}\supess_A\osc_{A,B_R(x)}X(A)<\infty$ and $\3X\3_{\osc,R}\le1$.
For all $t\in\R$, we apply ($\osc$-LSI) to the $\sigma(A)$-measurable random variable $e^{tX(A)/2}$,
\begin{align}\label{eq:LSI-conc-0bis}
\entm{e^{tX(A)}}\le C\,\expec{\int_{\R^d}\Big(\oscd{A,B_{R}(x)}e^{tX(A)/2}\Big)^2dx}.
\end{align}
By the inequality $|e^a-e^b|\le(e^a+e^b)|a-b|$ for all $a,b\in\R$, the integrand turns into
\begin{align}\nonumber
\Big(\oscd{A,S} e^{tX(A)/2}\Big)^2&\le2t^2\sup_{A,S}e^{tX(A)}\Big(\oscd{A,S}X(A)\Big)^2
\\
&\le {2t^2e^{tX(A)}}\exp\Big(t\oscd{A,S}X(A)\Big)\Big(\oscd{A,S}X(A)\Big)^2.\label{eq:bound-osc-exp}
\end{align}
Inserting this inequality into~\eqref{eq:LSI-conc-0bis} and using the assumptions on $X(A)$, we obtain
\begin{align*}
\entm{e^{tX(A)}}\le {Ct^2e^{tL}}\expecm{e^{tX(A)}}.
\end{align*}
Compared to the standard Herbst argument, we have to deal here with the additional exponential factor~$e^{tL}$. 
We may then appeal to~\cite[Corollary~2.12]{Ledoux-97} which indeed yields the desired Poisson concentration. 
We include a proof for the reader's convenience.
In terms of the Laplace transform $H(t)=\expecm{e^{tX(A)}}$, the above takes the form
\[tH'(t)-H(t)\log H(t)\le Ct^2e^{tL}H(t),\]
or equivalently,
\[\frac{d}{dt}\Big(\frac1t\log H(t)\Big) \le Ce^{tL},\]
and hence by integration
\[H(t) \le \exp\Big(\frac{Ct}L(e^{tL}-1)+t\frac{H'(0)}{H(0)}\Big) = e^{\frac{Ct}L(e^{tL}-1)+t\expec{X(A)}}.\]
The Markov inequality then implies for all $r,t\ge0$,
\begin{align}
\pr{X(A)\ge\expec{X(A)}+r}=\,\pr{e^{tX(A)}\ge e^{t\expec{X(A)}+tr}}&\le e^{-t\expec{X(A)}-tr}\expecm{e^{tX(A)}}\nonumber
\\
&\le\, e^{\frac{Ct}L(e^{tL}-1)-tr}.\label{eq:concentr-poiss-0}
\end{align}
Let $r\ge0$ be momentarily fixed, and denote by $t_*\ge0$ the value of $t\ge0$ that minimizes $f_r(t):=\frac{Ct}L(e^{tL}-1)-tr$, that is the (unique) solution $t_*\ge0$ of the equation 
\begin{equation}\label{e.add2}
Ce^{t_*L}=(Lr+C)/(1+t_*L)
\end{equation}
(note that $f_r$ is strictly convex, $f_r(0)=0$, and $f_r'(0)\le0$).
We now give two estimates on $f_r(t_*)$ depending on the value of $r$.
Assume first that $r\ge \frac{2e C}{L}$.
We may then compute
\[f_r(t_*):=\frac{Ct_*}L(e^{t_*L}-1)-t_*r\stackrel{\eqref{e.add2}}{=}-\frac{t_*^2(Lr+C)}{1+t_*L}.\]
Using the bound $2t_*L\ge t_* L+\log(1+t_*L)\stackrel{\eqref{e.add2}}{=}\log(1+Lr/C)$, 
and the fact that $t\mapsto -\frac{t^2(Lr+C)}{1+tL}$ is decreasing on $\R_+$, 
we obtain
\begin{align*}
f_r(t_*)\le-\frac{Lr+C}{2L^2}\frac{\log(1+Lr/C)^2}{2+\log(1+Lr/C)}.
\end{align*}
Hence, for $r\ge\frac{2eC}L$, we obtain using in addition $\log(1+Lr/C)\ge\log(1+2e)>9/5$,
\begin{align}\label{eq:concentr-poiss-pre}
f_r(t_*)\le-\frac{r}{5L}{\log\Big(1+\frac{Lr}{C}\Big)}.
\end{align}
We now turn to the case $0\le r\le \frac{2eC}{L}$. Comparing the minimal value $f_r(t_*)$ to the choice $t=\frac{r}{2eC}$, and using the bound $e^a-1\le ea$ for $a\in[0,1]$, we obtain for all $r\le\frac{2eC}{L}$,
\[f_r(t_*)\le f_r\Big(\frac{r}{2eC}\Big)= \frac {r}{2eL}\big(e^{\frac{rL}{2eC}}-1\big)-\frac{r^2}{2eC}\le-\frac {r^2}{4eC},\]
which yields, using that $\log(1+a)\le a$ for all $a\ge0$,
\[f_r(t_*)\le -\frac r{4eL}\log\Big(1+\frac{Lr}C\Big)\le -\frac r{11L}\log\Big(1+\frac{Lr}C\Big).\]
Combining this with~\eqref{eq:concentr-poiss-0} and~\eqref{eq:concentr-poiss-pre}, we conclude
\begin{align*}
\pr{X(A)\ge\expec{X(A)}+r}\le e^{-\frac{r}{11L}\log(1+\frac{Lr}C)},
\end{align*}
and the corresponding integrability result follows by integration.

\medskip
\step2 Proof of~(ii).

\nopagebreak
Let $A$ satisfy~($\osc$-MSG) with weight $\pi$, and let the random variable $X(A)$ satisfy~\eqref{e.conc-osc} for some $\kappa>0$, $p_0,\alpha\ge0$. Proposition~\ref{prop:pwsg}(ii) then yields for all $p\ge p_0$,
\begin{align*}
\expec{\big(X(A)-\expec{X(A)}\big)^{2p}}\le C^p p^{(2+\alpha) p}\kappa,
\end{align*}
or alternatively, for all $p\ge(2+\alpha)p_0$,
\begin{align*}
\expec{\Big(\big|X(A)-\expec{X(A)}\big|^{\frac{2}{2+\alpha}}\Big)^{p}}\le C^p p!\,\kappa.
\end{align*}
Summing this estimate over $p$, we obtain
\begin{align*}
\expec{\tilde\psi_{p_0,\alpha}\Big(\frac1C|X(A)-\expec{X(A)}|^{\frac2{2+\alpha}}\Big)}\,\le\,\kappa,
\end{align*}
where we have set $\tilde\psi_{p_0,\alpha}(u):=\sum_{n=0}^\infty \frac{u^{n+(2+\alpha)p_0}}{(n+(2+\alpha)p_0)!}$. Noting that $\tilde\psi_{p_0,\alpha}(u)\simeq_{p_0,\alpha}(1\wedge u)^{(2+\alpha)p_0}e^{u}$ holds for all $u\ge0$, the conclusion follows.
\qed


\section{Application to spatial averages}\label{chap:appl}

 This section is devoted to the proof of the results stated in Section~\ref{subsec:average}, that are, Propositions~\ref{prop:scl-averages}, \ref{prop:app-concentration}, and~\ref{prop:app-concentration-mixing}.

\subsection{Proof of Proposition~\ref{prop:scl-averages}: Scaling of spatial averages}
Let $L>0$. Given $\tilde \partial=\parfct{}$ or $\osc$, assumption~\eqref{eq:as-G} yields
\begin{multline*}
|\tilde \partial_{A,B_{\ell+1}(x)} X_L| \,\lesssim\, L^{-d}\int_{\R^d}e^{-\frac1L|y|}e^{-\frac1C(|x-y|-\ell)_+}dy \,\lesssim\, L^{-d}\int_{\R^d}e^{-\frac1L|y|}e^{-\frac1{C(\ell+1)}|x-y|}dy\\
\lesssim\, L^{-d}(L\wedge(\ell+1))^d e^{-\frac1{C(L+\ell+1)}|x|},
\end{multline*}
so that  the multiscale Poincaré inequality yields
\begin{align*}
\var{X_L} \,&\lesssim\, \int_0^\infty  \int_{\R^d} L^{-2d}(L\wedge(\ell+1))^{2d}e^{-\frac1{C(L+\ell+1)}|x|}dx\,(\ell+1)^{-d} \pi(\ell)d\ell\\
&\lesssim\, \int_0^\infty  L^{-2d}(L\wedge(\ell+1))^{2d}(L+\ell)^d(\ell+1)^{-d} \pi(\ell)d\ell\\
&\lesssim\, L^{-d}\int_0^L(\ell+1)^{d}\pi(\ell)d\ell+\int_L^\infty \pi(\ell)d\ell.
\end{align*}
An integration by parts yields $\pi_*(L)^{-1}\simeq L^{-d} \int_0^L  \pi(\ell) \ell^{d}d\ell+\int_L^\infty\pi(\ell)d\ell$, and the conclusion follows.
\qed

\subsection{Proof of Proposition~\ref{prop:app-concentration}: Linear concentration}

We split the proof into three steps.
We start with the proofs of  \eqref{e.app-conc1},~\eqref{e.app-conc2}, and~\eqref{e.app-conc3-0}, which directly follow from Propositions~\ref{prop:conc} and~\ref{prop:conc-osc}(ii).
The proof of estimates~\eqref{e.app-conc3} and~\eqref{e.app-conc4} is more subtle and is based on a fine tuning of the Herbst argument using specific features of the random variable~$Z_L$.

\medskip

\step1 Proof of \eqref{e.app-conc1},~\eqref{e.app-conc2}, and~\eqref{e.app-conc3-0}.

If $A$ satisfies ($\parfct{}$-MSG), we let the Lipschitz norm $\3\cdot\3_{\parfct{},\pi}$ be defined as in the statement of Proposition~\ref{prop:conc}. The same computation as in the proof of Proposition~\ref{prop:scl-averages} ensures that the random variable $Z_L:=\pi_*(L)^{1/2}X_L=\pi_*(L)^{1/2}\fint_{Q_L}(F-\expec{F})$ satisfies
\begin{align*}
\3Z_L\3_{\parfct{},\pi}\,\lesssim\, 1.
\end{align*}
Hence, estimates~\eqref{e.app-conc1} and~\eqref{e.app-conc2} follow from Proposition~\ref{prop:conc}.
We now turn to the proof of~\eqref{e.app-conc3-0}.
If $A$ satisfies ($\osc$-MSG) with weight $\pi(\ell)\lesssim(\ell+1)^{-\beta-1}$, $\beta>0$, we compute for all $p\ge p_0>\frac\beta d$, using assumption~\eqref{eq:as-G},
\begin{multline*}
\expec{\int_0^\infty\bigg(\int_{\R^d}\Big(\oscd{A,B_\ell(x)}X_L\Big)^2dx\bigg)^p(\ell+1)^{-dp-\beta-1}d\ell}\\
\lesssim L^{-2dp}\int_1^\infty(L+\ell)^{dp}(L\wedge \ell)^{2dp}\ell^{-dp-\beta-1}d\ell ~\lesssim~\big(1+(dp_0-\beta)^{-1}\big)\, L^{-\beta}.
\end{multline*}
Then applying Proposition~\ref{prop:conc-osc}(ii) and optimizing the choice of $p_0>\frac\beta d$, the result~\eqref{e.app-conc3-0} follows.

\medskip

\step2 Proof of~\eqref{e.app-conc3}.

\nopagebreak
Let $L\ge1$, and define $Z_L:=L^{d/2} X_L$.
As in the proof of Proposition~\ref{prop:scl-averages}, assumption~\eqref{eq:as-G} yields
\begin{equation}
\begin{array}{c}
\osc_{A,B_\ell(x)} Z_L \,\lesssim\, L^{-\frac d2} (L\wedge(\ell+1))^de^{-\frac1{C(L+\ell+1)}|x|}.
\end{array}
\label{eq:compute-ZL-osc}
\end{equation}
We make use of a variant of the Herbst argument as in~\cite[Section~4]{Bobkov-Ledoux-97} (see also~\cite[Section~2.5]{Ledoux-97}).
For all $t\ge0$ we apply ($\osc$-MSG) to the random variable $\exp(\frac12tZ_L)$: using the inequality $|e^a-e^b|\le(e^a+e^b)|a-b|$ for all $a,b\in\R$, we obtain
\begin{align*}
\varm{e^{\frac12tZ_L}}&\le\int_0^\infty\int_{\R^d}\expec{\Big(\oscd{A,B_\ell(x)}e^{\frac12tZ_L}\Big)^2}dx\,(\ell+1)^{-d}\pi(\ell)d\ell\\
&\lesssim t^2\,\expec{e^{tZ_L}}\supessd{A}\int_0^\infty\int_{\R^d} e^{t\oscd{A,B_\ell(x)}Z_L}\Big(\oscd{A,B_\ell(x)}Z_L\Big)^2dx\,e^{-\frac1C\ell^\beta}d\ell,
\end{align*}
and hence, in terms of the Laplace transform $H_L(t):=\expec{e^{tZ_L}}$,
\begin{equation*}
{H_L(t)-H_L(t/2)^2\le ~ t^2H_L(t) \supessd{A} \int_0^\infty \int_{\R^d}  e^{t \osc_{A,B_{\ell}(x)} Z_L} \Big(\osc_{A,B_\ell(x)}Z_L\Big)^2 dx\, e^{-\frac1C \ell^\beta} d\ell}.
\end{equation*}
Using the property~\eqref{eq:compute-ZL-osc} of the random variable $Z_L$, we find
\begin{eqnarray}
\lefteqn{H_L(t)-H_L(t/2)^2}\nonumber\\
&\lesssim& t^2H_L(t)  \int_1^\infty \Big(\frac{L\wedge \ell}{\sqrt L}\Big)^{2d}\exp\bigg(Ct\Big(\frac{L\wedge \ell}{\sqrt L}\Big)^d-\frac{\ell^\beta}C\bigg)\int_{\R^d}e^{-\frac1{C(L+\ell+1)}|x|}dx\,d\ell\nonumber\\
&\lesssim& t^2H_L(t) \int_1^\infty (L+\ell)^d\,\Big(\frac{L\wedge \ell}{\sqrt L}\Big)^{2d}\exp\bigg(Ct\Big(\frac{L\wedge \ell}{\sqrt L}\Big)^d-\frac{\ell^\beta}C\bigg)\,d\ell\nonumber\\
&\lesssim& t^2H_L(t)\bigg(\int_0^L \exp\bigg(Ct \Big(\frac{\ell}{\sqrt{L}}\Big)^d-\frac{\ell^\beta}C\bigg)d\ell+L^d e^{CtL^\frac d2}\int_L^\infty e^{-\frac1C \ell^\beta} d\ell\bigg).\label{eq:est-HL-herbst-sg}
\end{eqnarray}
Without loss of generality we may assume that $\beta\le \frac d2$ (the statement~\eqref{e.app-conc3} is indeed not improved for $\beta>\frac d2$).
We then restrict to
\begin{align}\label{eq:choice-T-sg}
0\le t\le T:=\frac1{K} L^{\beta-\frac d2},
\end{align}
for some $K\gg 1$ to be chosen later (with in particular $K\ge2C^2$). As a consequence of $\beta\le\frac d2$, this choice yields $T\le K^{-1}$.
On the one hand, for all $0\le \ell \le L$ and all $0\le t\le T$, the choice of $T$ with $K\ge2C^2$ yields
$$
Ct \Big(\frac{\ell}{\sqrt{L}}\Big)^d-\frac{\ell^\beta}C\,=\,-\frac{L^\beta}C \bigg(\Big(\frac \ell L\Big)^\beta-\frac{C^2t}{L^{\beta-\frac d2}} \Big(\frac\ell L\Big)^d\bigg)\le\,-\frac{L^\beta}C \bigg(\Big(\frac \ell L\Big)^\beta-\frac12 \Big(\frac\ell L\Big)^d\bigg)\le -\frac{\ell^\beta}{2C},
$$
and hence
$$
\int_0^L \exp\bigg(Ct \Big(\frac{\ell}{\sqrt{L}}\Big)^d-\frac{\ell^\beta}C\bigg)d\ell\,\lesssim\, \int_0^\infty e^{ -\frac{\ell^\beta}{2C}}d\ell \,\lesssim\, 1.
$$
On the other hand, for all $0\le t\le T$, the choice of $T$ with $K\ge2C^2$ yields
$$
L^d e^{CtL^\frac d2}\int_L^\infty e^{-\frac1C \ell^\beta} d\ell\,\lesssim\, \exp\Big(CtL^\frac d2-\frac{L^\beta}{2C}\Big)\le \exp\Big(\frac{CL^\beta}K-\frac{L^\beta}{2C}\Big)\le1.
$$
Injecting these estimates into~\eqref{eq:est-HL-herbst-sg}, we obtain for all $0\le t\le T$,
\begin{align*}
H_L(t)-H_L(\tfrac t2)^2&\le Ct^2H_L(t),
\end{align*}
and hence
\begin{align*}
H_L(t)\le \frac{H_L(\tfrac t2)^2}{1-Ct^2}.
\end{align*}
Applying the same inequality for $t/2$, iterating, and noting that $H_L(2^{-n}t)^{2^n}\to e^{t\expec{Z_L}}=1$ as $n\uparrow\infty$, we obtain for all $0\le t\le T$,
\begin{align*}
H_L(t)\le \prod_{n=0}^{\infty}\big(1-C(2^{-n}t)^2\big)^{-2^{n}}.
\end{align*}
For $K$ large enough such that $CT^2\le CK^{-2}\le \frac12$, the inequality $\log(1-x)\ge -2x$ for all $0\le x\le\frac12$ then yields for all $0\le t\le T$,
\begin{align*}
\log H_L(t)\le -\sum_{n=0}^{\infty}2^n\log\big(1-C(2^{-n}t)^2\big)\le 2Ct^2\sum_{n=0}^{\infty}2^{-n}\lesssim t^2,
\end{align*}
and thus $H_L(T)\le e^{CT^2}$. Using Markov's inequality and the choice~\eqref{eq:choice-T-sg} of $T$, we deduce for all $r\ge0$,
\[\pr{Z_L>r}\le e^{-Tr+CT^2}=\exp\Big(-\frac{L^{\beta-\frac d2}r}K+\frac{C}{K^2}L^{2\beta -d}\Big).\]
With the choice $r=\delta L^\frac d2$ for $\delta>0$, this turns into
\[\pr{X_L>\delta}\le \exp\Big(-\frac\delta KL^\beta+\frac C{K^2}L^{2\beta-d}\Big)\le\exp\bigg(-\frac1K\Big(\delta-\frac CK\Big)L^\beta\bigg).\]
Choosing $K\simeq1\vee\delta^{-1}$ large enough, the desired estimate~\eqref{e.app-conc3} follows.

\medskip

\step3 Proof of~\eqref{e.app-conc4}.
\nopagebreak

Let $L\ge1$, and define $Z_L:=L^{d/2}X_L$. We make use of the classical Herbst argument as presented e.g.\@ in~\cite[Section~5.1]{Ledoux-01}.
For all $t\ge 0$ we apply ($\osc$-MLSI) to the random variable $\exp(\frac12 tZ_L)$,
\begin{align*}
\ent{e^{tZ_L}} \,&\le\, \int_0^\infty \int  \expec{\Big(\osc_{A,B_{\ell}(x)} e^{\frac12 tZ_L}\Big)^2}dx\,(\ell+1)^{-d}\pi(\ell)d\ell.
\end{align*}
Estimating the RHS as in~\eqref{eq:est-HL-herbst-sg}, we obtain in terms of $H_L(t):=\expecm{e^{tZ_L}}$,
\begin{align*}
\frac{d}{dt}\Big(\frac1t\log H_L(t)\Big)\,&\lesssim \, \int_0^L \exp\bigg(Ct \Big(\frac{\ell}{\sqrt{L}}\Big)^d-\frac{\ell^\beta}C\bigg)d\ell+L^d e^{CtL^\frac d2}\int_L^\infty e^{-\frac1C \ell^\beta} d\ell.
\end{align*}
Without loss of generality we may assume that $\beta\le d$ (the statement~\eqref{e.app-conc4} is indeed not improved for $\beta>d$).
We then restrict to
\begin{align}\label{eq:choice-T}
0\le t\le T:=\frac1{K} L^{\beta-\frac d2},
\end{align}
for some $K\gg 1$ to be chosen later (with in particular $K\ge2C$).
Arguing as in Step~1, we obtain for all $0\le t\le T$,
\[\frac{d}{dt}\Big(\frac1t\log H_L(t)\Big)\,\lesssim\, 1,\]
which yields by integration with respect to $t$ on $[0,T]$,
\begin{align*}
\frac1T \log H_L(T)=\frac1T \log H_L(T)-\expec{Z_L}\,&\lesssim \, T,
\end{align*}
that is, $H_L(T) \,\le \, e^{C T^2}$. The desired estimate~\eqref{e.app-conc4} then follows as in Step~1, using Markov's inequality and choosing $K$ large enough.
\qed

\subsection{Proof of Proposition~\ref{prop:app-concentration-mixing}: Linear concentration}

Without loss of generality we assume that $\supess_A |f(A)|\le 1$, which implies $\pr{|X_L|> 1}=0$. It is then sufficient to establish the result for $0<\delta\le1$.
We split the proof into two steps.
In the first step we prove the result in the case when the random variable $f(A)$ is exactly $1$-local.
We then extend the result in Step~2 when $f(A)$ is only approximately local in the sense~\eqref{eq:as-G0}.
Since (MCI) implies $\alpha$-mixing by Proposition~\ref{prop:sgmix}(iii), it is enough to prove the result under the sole assumption of $\alpha$-mixing.

\medskip

\step1 Exactly $1$-local random variable $f(A)$.
\nopagebreak

In this step we assume in addition that $f(A)$ is $\sigma(A|_{B_1})$-measurable, and we prove that for all $\delta,L>0$,
\begin{equation}\label{conc-exact-local}
\pr{X_L>\delta}\le C\exp\Big(-\frac{\delta^2}CL^{\frac{d\beta}{d+\beta}}\Big).
\end{equation}
Let $p\ge1$ be an integer and let $R>0$. Setting
\[E_{R,p}:=\{(x_1,\ldots,x_{2p})\in(\R^d)^{2p}:|x_1-x_j|>R,\,\forall j\ne1\},\]
and noting that for all $x$ the random variable $F(x)$ is $\sigma(A|_{B(x)})$-measurable,
$\alpha$-mixing leads to
\begin{multline}\label{eq:estimate-prod-far}
L^{-2dp}\bigg|\int\!\!\ldots\!\!\int_{E_{R,p}}e^{-\frac1L\sum_{i=1}^{2p}|x_i|}\,\expec{(F(x_1)-\expec{F})\ldots(F(x_{2p})-\expec{F})}dx_1\ldots dx_{2p}\bigg|\\
\le C^pL^{-2dp}\int\!\!\ldots\!\!\int_{E_{R,p}}e^{-\frac1L\sum_{i=1}^{2p}|x_i|}\,\tilde\alpha\bigg(R-2,2+2\sum_{i=1}^{2p}|x_i|;A\bigg)dx_1\ldots dx_{2p}\\
\le C^pL^{-2dp}e^{-\frac1CR^\beta}\int_{\R^{2dp}} e^{-\frac1L|x|}(1+|x|)^Cdx\le C^pL^Ce^{-\frac1CR^\beta}.
\end{multline}
Using this estimate, we compute
\begin{multline*}
\expecm{X_L^{2p}}=L^{-2dp}\int_{\R^d}\!\!\!\ldots\!\int_{\R^d}e^{-\frac1L\sum_{i=1}^{2p}|x_i|}\,\expec{(F(x_1)-\expec{F})\ldots(F(x_{2p})-\expec{F})}dx_1\ldots dx_{2p}\\
\le C^pL^{C}e^{-\frac1CR^\beta}+C^pL^{-2dp}\int_{\R^d}\!\!\!\ldots\!\int_{\R^d}e^{-\frac1L\sum_{i=1}^{2p}|x_i|}\,\mathds1_{\forall i,\,\exists j\ne i:\,|x_i-x_j|\le R}\,dx_1\ldots dx_{2p}.
\end{multline*}
We consider the partitions $P:=\{P_1,\ldots,P_{N_P}\}$ of the index set $[2p]:=\{1,\ldots,2p\}$ into nonempty subsets of cardinality $\ge2$ (that is, $\cup_jP_j=[2p]$, $\sharp P_j\ge2$ for all $j$, and $P_j\cap P_l=\varnothing$ for all $j\ne l$), and we use the notation $P\vdash_2[2p]$ for such partitions. The above then takes the form
\begin{align*}
\expecm{X_L^{2p}}&\le C^pL^{C}e^{-\frac1CR^\beta}+C^pL^{-2dp}\sum_{P\vdash_2[2p]}L^{dN_P}R^{d(2p-N_P)}.
\end{align*}
Since for all $1\le k\le p$ the number of partitions $P\vdash_2[2p]$ with $N_P=k$ is bounded by the Stirling number of the second kind $\lbrace{2p\atop k}\rbrace\le\frac12\binom{2p}{k}k^{2p-k}\le C^pp^{2p}k^{2(p-k)}(2p-k)^{-(2p-k)}$, we deduce
\begin{align*}
\expecm{X_L^{2p}}&\le C^pL^{C}e^{-\frac1CR^\beta}+C^p\Big(\frac RL\Big)^{dp}\sum_{k=1}^p\frac{p^{2p}k^{2(p-k)}}{(2p-k)^{2p-k}}\Big(\frac RL\Big)^{d(p-k)},
\end{align*}
and hence by Markov's inequality, for all $\delta>0$,
\begin{align}\label{eq:bound-proba-XL}
\pr{X_L>\delta}\le \delta^{-2p}C^pL^{C}e^{-\frac1CR^\beta}+\delta^{-2p}C^p\Big(\frac RL\Big)^{dp}\sum_{k=1}^p\frac{p^{2p}k^{2(p-k)}}{(2p-k)^{2p-k}}\Big(\frac RL\Big)^{d(p-k)}.
\end{align}
Recall that we may restrict to $0<\delta \le 1$.
Choosing $R=L^\alpha$, $p=\delta^2C_0^{-1}L^{\alpha\beta}$, and $\alpha=\frac{d}{d+\beta}$, for some $C_0\simeq1$ large enough, the estimate~\eqref{eq:bound-proba-XL} above leads to
\begin{align*}
\pr{X_L>\delta}\le Ce^{-\frac1CL^{\alpha\beta}}+\delta^{-2p}C^pL^{-\alpha\beta p}\sum_{k=1}^p\frac{p^{p+k}k^{2(p-k)}}{(2p-k)^{2p-k}}.
\end{align*}
Noting that the summand is increasing in $k$, and using the choice of $p$ with $C_0$ large enough, we deduce 
\begin{align}
\pr{X_L>\delta}&\le Ce^{-\frac1CL^{\alpha\beta}}+\delta^{-2p}C^pL^{-\alpha \beta p}p^p\le Ce^{-\frac1C\delta^2L^{\alpha\beta}},
\end{align}
from which the desired result~\eqref{conc-exact-local} follows.

\medskip

\step2 Approximately $1$-local random variable $f(A)$.

For all $r>1$, we define the ($r$-local) random variable $f_r(A):=\expeC{f(A)}{A|_{B_r}}$, and we set $F_r(x):=f_r(A(\cdot+x))$ and $X_{r,L}:=L^{-d}\int_{\R^d}e^{-\frac1L|y|}(F_r(y)-\expec{F_r})dy$.
The approximate locality assumption~\eqref{eq:as-G0} implies a.s.\@ for all $r,L>0$,
\begin{align}\label{eq:as-G0-appl}
|X_{r,L}-X_L|\le Ce^{-\frac rC}.
\end{align}
Setting $\tilde F_r(x):=F(rx)$ and $A_r(x):=A(rx)$, we note that for all $x\in\R^d$ the random variable $\tilde F_r(x)$ is $\sigma(A|_{B_r(rx)})$-measurable, that is, $\sigma(A_r|_{B(x)})$-measurable. For all $r\ge1$, the $\alpha$-mixing assumption on $A$ implies that the contracted random field $A_r$ satisfies $\alpha$-mixing with coefficient 
$$
\tilde\alpha_r(\ell,D;A_r):=\Big((1+r D)^C\exp(-\frac1C(r\ell)^\beta)\Big) \wedge 1 \,\le \, C(1+r D)^C\exp(-\frac1C(r\ell)^\beta),
$$
so that the $\alpha$-mixing coefficient for $r\ge 1$ can basically be estimated by the one for $r=1$. 
We may therefore apply Step~1 in the following form for all $\delta,L>0$ and all $r\ge1$,
\[\pr{X_{r,L}>\delta}=\pr{\fint_{Q_{L/r}}(\tilde F_r-\expecm{\tilde F_r})>\delta}\le C\exp\bigg(-\frac{\delta^2}C\Big(\frac Lr\Big)^{\frac{d\beta}{d+\beta}}\bigg),\]
where the constant $C\ge1$ is independent of $r$.
Combining this with~\eqref{eq:as-G0-appl} and choosing $r:=C|\log(\frac{\delta}{eC})|\ge1$, we obtain for all $0<\delta\le 1$ and $L>0$,
\begin{align*}
\pr{X_{L}>\delta}\le\pr{X_{r,L}>\delta-Ce^{-\frac rC}}\le\pr{X_{r,L}>\tfrac\delta2}\le C\exp\bigg(-\frac{\delta^2}{C}\,\Big(\frac{L}{|\log(\frac{\delta}{eC})|}\Big)^{\frac{d\beta}{d+\beta}}\bigg),
\end{align*}
and the conclusion follows.
\qed


\section*{Acknowledgements}
The work of MD is supported by F.R.S.-FNRS (Belgian National Fund for Scientific Research) through a Research Fellowship.
The authors acknowledge financial support from the European Research Council under
the European Community's Seventh Framework Programme (FP7/2014-2019 Grant Agreement
QUANTHOM 335410).

\bigskip
\bibliographystyle{plain}

\def\cprime{$'$} \def\cprime{$'$}

\end{document}